\documentclass[conference]{IEEEtran}
\usepackage{cite}

\ifCLASSINFOpdf
   \usepackage[pdftex]{graphicx}
   \DeclareGraphicsExtensions{.pdf,.jpeg,.png}
\else
\fi
\usepackage{amssymb}
\usepackage[cmex10]{amsmath}

\usepackage{diagbox}
\usepackage[flushleft]{threeparttable}
\usepackage{soul}
\usepackage{color}
\usepackage{optidef}
\usepackage{enumitem}

\begin{document}
%
\title{A Framework to Utilize DERs' VAR Resources to Support the Grid in an Integrated T-D System}

\author{\IEEEauthorblockN{Ankit Singhal}
\IEEEauthorblockA{Department of Electrical and\\Computer Engineering\\
Iowa State University\\
Ames, Iowa 50010\\
Email: ankit@iastate.edu}
\and
\IEEEauthorblockN{Venkataramana Ajjarapu}
\IEEEauthorblockA{Department of Electrical and\\Computer Engineering\\
Iowa State University\\
Ames, Iowa 50010\\
Email: vajjarap@iastate.edu}}


%


\maketitle

\begin{abstract}
 Increasing penetration of inverter-based distributed energy resources (DERs) opens up interesting opportunities for the transmission systems. We present a hypothesis that the numerous DERs in var control mode can be seen as geographically distributed var devices (mini-SVCs) and if controlled properly, can be exploited to increase system flexibility by providing local var support to the grid as an ancillary service. Based on this premise, a var support framework is proposed in this paper. It utilizes a novel D-OPF formulation for unbalanced three-phase feeders enabling the estimation of the maximum var support that can be provided by the DERs to the grid at different operating points without 
compromising the distribution network performance. Further, a co-simulation method is developed to investigate the true impact of the proposed DER var support on the grid in an integrated Transmission-Distribution (T-D) system.

\end{abstract}
\IEEEpeerreviewmaketitle

\section{Introduction}
The conducive environment for distributed energy resources (DER) growth is pushing its penetration to as high as 100\% \cite{bank_high_2013}. 
Since the distribution system is not physically isolated from the bulk transmission systems, the cumulative influence of the increasing DERs on the transmission grid performance should not be ignored. 
A discussion on this topic has started relatively recently \cite{perez-arriaga_transmission_2016} and there is a growing need to assess DERs' impacts on the grid in the literature. 

While the DER integration can impact the grid both adversely and positively, the purpose of this paper is to explore the opportunities offered by the DERs to improve transmission system performance. In particular, we present a hypothesis that, from the grid perspective, thousands of DER devices with volt/var control capability can be seen as the geographically distributed var resources (\textit{`mini- static voltage compensators}) to improve system performance. It is known that the bulk transmission system needs to install its own reactive power devices such as static voltage compensators (SVC) and capacitors at certain locations to enhance the system flexibility and voltage stability. Usually, these devices are very costly whereas due to the distributed nature, DERs can provide more flexible and localized var support to the bulk system at less economic cost, if controlled 
properly. Based on this premise, 
we propose a DER var support framework that utilizes the DERs volt/var capability for the benefits of the grid and verifies its impact in an integrated transmission-distribution (T-D) network.

Mainly, the DERs can support the grid by influencing the net real and reactive power flows at the feeder substation through the var injection at their local node of connection; however, this can adversely affect the distribution feeder voltages. Note that all the utilities are enforced to maintain feeder voltages within the allowable range by ANSI standard \cite{noauthor_ansi_2016}. Therefore, 
we propose a distributed optimal power flow (D-OPF) that estimates the maximum var support that a feeder can provide at the substation at different operating points throughout the day without violating its own operational limits. 
There have been few studies to utilize DERs to improve distribution system performance at the substation; for instance, \cite{dallanese_optimal_2017} minimizes the real power demand of the feeder using DERs and \cite{keane_enhanced_2011} minimize the var demand using fixed power factor mode of DERs. However, these studies do not consider unbalanced distribution feeders and impact of DER on the grid. 

One of the major reasons of lack of studies on the DERs impact on the grid is lack of an appropriate co-simulation platform which can solve the integrated transmission and distribution system, though there have been some recent attempts to develop T-D co-simulation methods and platforms \cite{sun_master_2015,hansen_bus.py:_2015,palmintier_igms:_2017}. The literature in co-simulation is still in a nascent phase and different platforms are being developed based on different open source solvers and the application of interest. In the proposed framework, we utilize a T-D co-simulation platform based on the widely accepted power flow solvers (Matpower and GridlabD) and couple it with D-OPF to study and verify the true impact of DERs var support on the grid. The integrated system allows to include and observe the distribution system details as well as the changing substation voltage behavior which are not possible in traditional aggregated load modeling.

Overall, the objective of this work is to investigate the var support opportunities provided by the DERs for the grid benefits by proposing a DER var support framework. The support framework has two main novel aspects i.e. D-OPF and the T-D co-simulation which intend to provide following unique contributions: 1) To estimate the \textit{`day ahead maximum var support curve'} for the grid while enforcing distribution system operating limits; 2) To verify the true impact of var support on the grid
 using an integrated T-D co-simulation; 3) To confirm our proposition that the DERs can be exploited as mini-SVCs to provide flexibility and ancillary services to the grid; and 4) To provide a general framework which enables further investigation of the DERs utilization for various grid support applications in an integrated T-D environment.

\section{Preliminaries and System Modeling}
Consider a radial distribution network with $N+1$ nodes represented by a tree graph $\mathcal{T}=(\mathcal{N},\mathcal{E}$), where $\mathcal{N}:=\{0,1,\cdots,N\}$ is a set of distribution nodes and $\mathcal{E}:=\{(i,j)\}$ is a set of line segments. The subset $\mathcal{N}_j$ is a collection of all immediate downstream neighboring buses of node $j$. The feeder head or the substation is denoted by node 0 which is assumed to be the reference voltage for the feeder.
The basic KVL and KCl equations for for buses $j$ and $k$ in an unbalanced distribution system can be written as \cite{kersting_distribution_2017}:
 \setlength\abovedisplayskip{3pt}
 \setlength\belowdisplayskip{3pt}
\begin{equation}
\label{eq:kvl}
 \mathbb{V}_j = \mathbb{V}_k + \mathbb{Z}_{jk}\mathbb{I}_k
\end{equation}
\begin{equation}
\label{eq:kcl}
 \mathbb{I}_j = \mathbb{I}^{inj}_j + \sum_{k\in \mathcal{N}_j } \mathbb{I}_k
\end{equation}
\vspace{-1mm}
Where, $\mathbb{V}_j=[V_a V_b V_c]_j^T$ represent the vector of voltage phasors at node $j$ and $\mathbb{I}_j=[I_a I_b I_c]_j^T$ represent the vector of current phasors entering at node $j$. Similarly, the net current phasor injected at node $j$ is denoted by $\mathbb{I}^{inj}_j$. 
$\mathbb{Z}_{jk}$ represents the 3 phase impedance matrix ($ 3\times 3$) of the line segment $jk$ as shown in \cite{kersting_distribution_2017}. Let's define the complex power entering at phase $\phi$ of node $j$ as $S_{\phi,j}=V_{\phi,j}I_{\phi,j}^*$ and complex net power injected at phase $\phi$ of node $j$ as $s_{\phi,j}=V_{\phi,j}I_{\phi,j}^{inj*}$ for all three phases $\phi\in\{a,b,c\}$. Now, by following the \textit{LinDistFlow} formulation as described in \cite{baran_optimal_1989,arnold_optimal_2016}, (\ref{eq:kvl})-(\ref{eq:kcl}) can be re-written as,
\begin{equation}
\label{eq:LinDist1a}
\mathbb{S}_j\approx -s_j+\sum_{k \in \mathcal{N}_j}\mathbb{S}_k
\end{equation}
\begin{equation}
\label{eq:LinDist1b}
\mathbb{V}_j\mathbb{V}_j^* \approx \mathbb{V}_k\mathbb{V}_k^* + 2 \textbf{Re} \{\mathbb{V}_k[S_a V_a^{-1} S_b V_b^{-1} S_c V_c^{-1}]_k\mathbb{Z}_{jk}^*\}
\end{equation}
Where $\mathbb{S}_j$ and $s_j$ are the vectors of complex phasors $S_{\phi,j}$ and $s_{\phi,j}$, respectively.  Here line losses are neglected which introduce a relatively small error in the modeling as indicated by \cite{farivar_equilibrium_2013}. However, to increase accuracy, we will be adding a constant loss term $L_j$ in (\ref{eq:LinDist1a}). To further linearize the system, another approximation of constant ratio of voltage phasors is assumed \cite{arnold_optimal_2016}. Now, we define the vector of squared of voltage magnitude as a new variable $\mathbb{Y}_j\! =\! \mathbb{V}_j\mathbb{V}_j\!=\![y_a \ y_b \ y_c]$, 
 which will be used in the formulation to keep it linear rather than $\mathbb{V}_j$. Expanding the impedance matrix entries into resistance and reactances as $Z_{\phi \psi}=r_{\phi \psi}+jx_{\phi \psi}$, and complex line flows into real and reactive power flows as $S_{\phi,j}=P_{\phi,j}+jQ_{\phi,j}$, (\ref{eq:LinDist1b}) can be written as following linear matrix equation:
 \begin{equation}
  \label{eq:LinDist2a}
 \mathbb{Y}_j = \mathbb{Y}_k - \mathbb{M}_{jk}^p\mathbb{P}_k - \mathbb{M}_{jk}^Q\mathbb{Q}_k 
 \end{equation}
 
 
Where, $\mathbb{P}_j\!=\![P_a \ P_b \ P_c]$, 
$\mathbb{Q}_j\!=\! [Q_a \ Q_b \ Q_c]$, and  $M_{jk}^P$ and $M_{jk}^Q$ can be obtained from \cite{arnold_optimal_2016}.
  
 (\ref{eq:LinDist1a}) is re-written with a constant loss term as:
 \begin{equation}
  \label{eq:LinDist2d}
\mathbb{S}_j\approx -s_j+\sum_{k \in \mathcal{N}_j}\mathbb{S}_k + L_j
 \end{equation}
The linearized model of an unbalanced 3 phase distribution system is represented by the equations (\ref{eq:LinDist2a})-(\ref{eq:LinDist2d}). The loss term $L_j$ can be estimated based on the offline study of the base operating point as indicated in \cite{zhu_fast_2016}. 

Let's assume the DERs are located at the nodes collected in a subset $\mathcal{G}\subseteq \mathcal{N}$. In this case, only inverter based DERs are considered such as solar PV. The net power injection of real and reactive power at each node $j \in \mathcal{G}$ is denoted by $p_{\phi,j}=p_{\phi,j}^G - p^L_{\phi,j}$ and $q_{\phi,j}=q_{\phi,j}^{inv} - q^L_{\phi,j}$ respectively. Superscript $G$ and $L$ denote generation and consumption of power respectively, whereas, $q^{inv}_{\phi,j}$ denotes the controllable reactive power injection by the inverter at phase $\phi$ of the local node $j$. For all other $j\in \mathcal{N-G}$, $p_{\phi,j}^G$ and $q_{\phi,j}^{inv}$ are considered zero. 
Only constant power loads are taken here and the capacitors are modeled as reactive power loads.  

\section{Proposed VAR Support Framework}
\begin{figure}[b]
	\centering
    \vspace{-4mm}
	\includegraphics[trim=-0.5in 0in 0in 0in,width=3.1in]{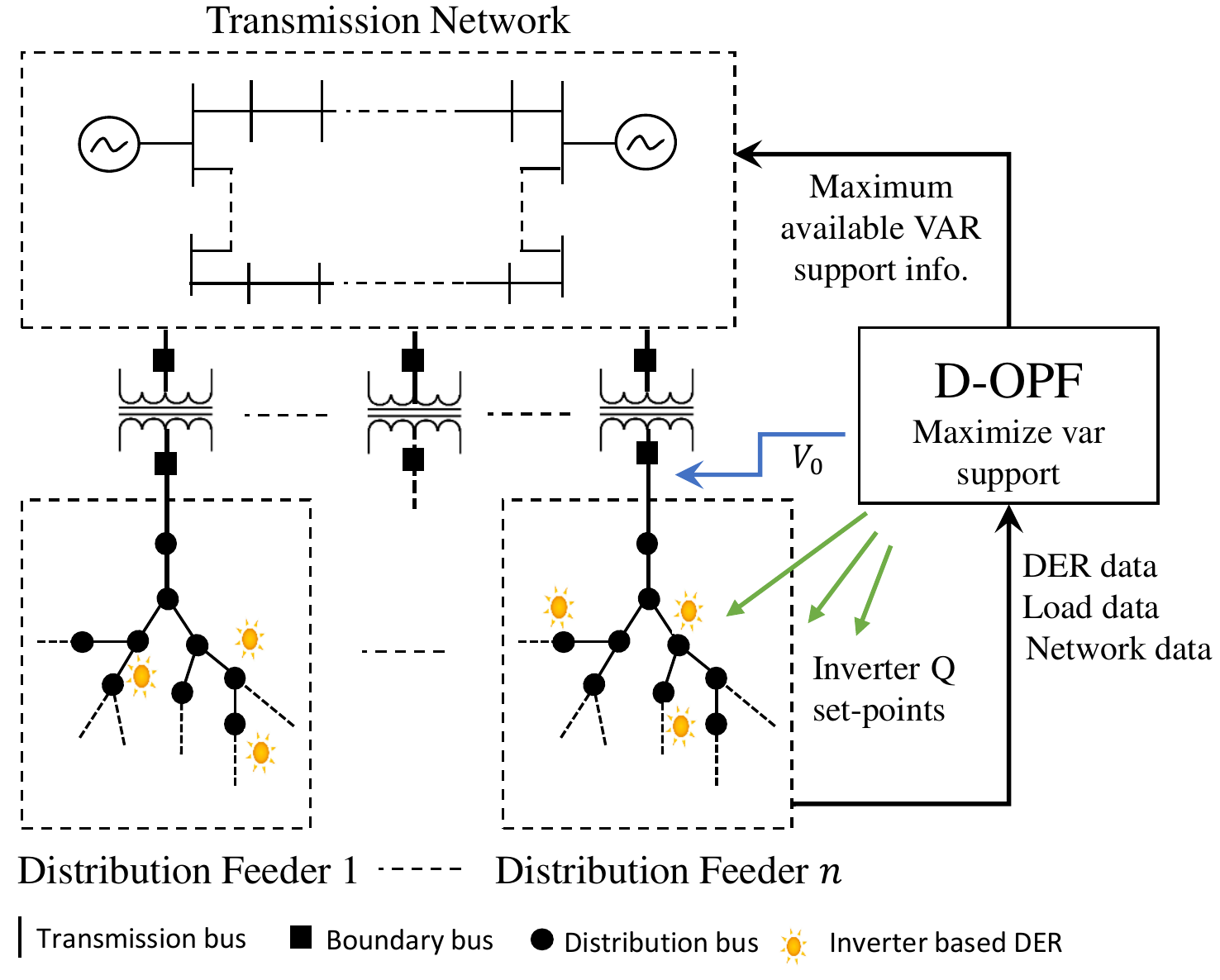}
    \vspace{-3mm}
    \caption {The DER var support framework for an integrated T-D system }
    \label{fig:overall}
\end{figure}

 \figurename \ref{fig:overall} depicts the overall framework of providing var support to the grid from DERs in an integrated T-D system, proposed in this work. 
 The objective of this framework is to estimate the \textit{maximum var support capacity curve} and send it to the transmission system operator (TSO).  This is achieved through a D-OPF module which estimates the maximum available var support at the feeder substation at different operating points throughout the day while ensuring the distribution system is within its operational limits. Then the distribution system operator (DSO) collect the maximum available var support from all the connected feeders at different operating points and aggregates it in form of maximum var support capacity curve and sends it to the TSO so that the transmission controller can take the appropriate decision. Here we assume that the TSO has its own monitoring and control methods to request var support from the DSO in case of emergency. Another function of D-OPF module is to dispatch optimal inverter var set-points to individual DER devices in order to meet the var support requested by the grid; however, in this work we do not provide details of this functionality due to space limitations and only focus on developing a general framework with maximum var support estimation functionality. Other functions of the framework will be explored in the future studies. 
 
 \subsection{D-OPF Formulation}
 The objective of maximizing the available local var support from DERs is same as minimizing the net reactive power demand at the substation for the grid i.e. $Q^{net}_0$. The expression for $Q^{net}_0$ can be approximated as,
\begin{equation}
\small
\label{eq:q_net}
 Q^{net}_0 =\mathcal{L}(y_{\phi,j})+ \sum_{\phi \in \{a,b,c\}} \Big(\sum_{j\in \mathcal{G}}q_{\phi,j}^{inv}-\sum_{j\in \mathcal{N}}q_{\phi,j}^{L}\Big)
\end{equation}
Where, the first term $\mathcal{L}(\mathbb{Y}_j)$ represents total reactive power losses in the system which can be written as $\sum_j S_j^2/V_j^2.x_j$ for a balanced single phase system \cite{baran_network_1989}, however, for an unbalanced system, the loss expression becomes complicated due to interaction all the phases. To simplify it, an approximation is considered that assumes the influence of non-diagonal entries of $\mathbb{Z}$ negligible compared to the influence of diagonal entries while estimating the losses. Based on this assumption, $\mathcal{L}(\mathbb{Y}_j)$ can be written as:
 \begin{equation}
 \small
 \label{eq:loss}
\mathcal{L}(y_{\phi,j}) = \sum_{\phi \in \{a,b,c\}} \sum_{j\in \mathcal{N}} \frac{(P_{\phi,j}^2+Q_{\phi,j}^2)}{y_{\phi,j}} x_{\phi\phi,j}
\end{equation}
The second term in (\ref{eq:q_net}) represents the total net injection of var at each node due to loads, capacitors and DER inverters. Based on the already defined preliminaries, We define the D-OPF as:
\vspace{-2mm}
\begin{mini!}|l|[2]
		{\small q_{\phi,j}^{inv},y_{\phi,0}}{Q^{net}_0(y_{\phi,j},q_{\phi,j}^{inv})}
		{\label{opt}}{}
        \vspace{-1mm}
        \addConstraint{(\ref{eq:LinDist2a})-(\ref{eq:LinDist2d})}{,}{\forall j\in \mathcal{N}
        }
        \addConstraint{p_{\phi,j}=p_{\phi,j}^G - p^L_{\phi,j}}{,\label{opt:DGP_const}\quad}{\forall j\in \mathcal{N}}
        \addConstraint{q_{\phi,j}=q_{\phi,j}^{inv} - q^L_{\phi,j}}{,\label{opt:DGQ_const}\quad}{\forall j\in \mathcal{N}}
		\addConstraint{\underline{y}\leq y_{\phi,j}\leq \overline{y}}{,\label{opt:volt_const}\quad}{\forall j\in \mathcal{N}}
        \addConstraint{|q_{\phi,j}^{inv}|\leq \overline{q_{\phi,j}} }{,\label{opt:q_const}\quad}{\forall j\in \mathcal{G}}
\end{mini!}
Constraint (\ref{opt:volt_const}) ensures the voltages are within the ANSI limits \cite{noauthor_ansi_2016}. $\overline{y}$ and $\underline{y}$ are upper and lower allowable voltage limits, and are usually taken as $1.05^2$ and $0.95^2$, respectively. Constraint (\ref{opt:q_const}) manifest the hardware capacity limit of an inverter. The maximum available inverter var at any time instant ($\overline{q_{\phi,j}}$) depends on the inverter capacity, {\small $S^{inv}_j$}, and the real power generation at that time instant i.e. {\small $\overline{q_{\phi,j}}=\sqrt{S^{inv^2}_j-p_{\phi,j}^{G^2}}$}.
Although, DER real power generation curtailment is not advisable in normal situations, this framework allows the option of DER curtailment which provides more flexibility to the operator as discussed in the next section. The solution of the optimization (\ref{opt}) provides the optimal var set-point for each dispatchable DER inverter in form of $q_{\phi,j}^{inv^*}$ and optimal secondary side voltage set-point ($y_0$) in order to maximize the available var support to the grid. Note that the on-load tap changer (OLTC) can control the secondary side voltage through tap changes within a range which is implemented here. We assume the following LTC logic: each tap provides $\pm$0.01 pu regulation with maximum $\pm10$  taps with 0.01 pu as bandwidth.
\vspace{-1mm}
\subsection{Co-simulation Framework}
\begin{figure}
	\centering
	\includegraphics[trim=-0in 0in 0in 0in,width=2.5in]{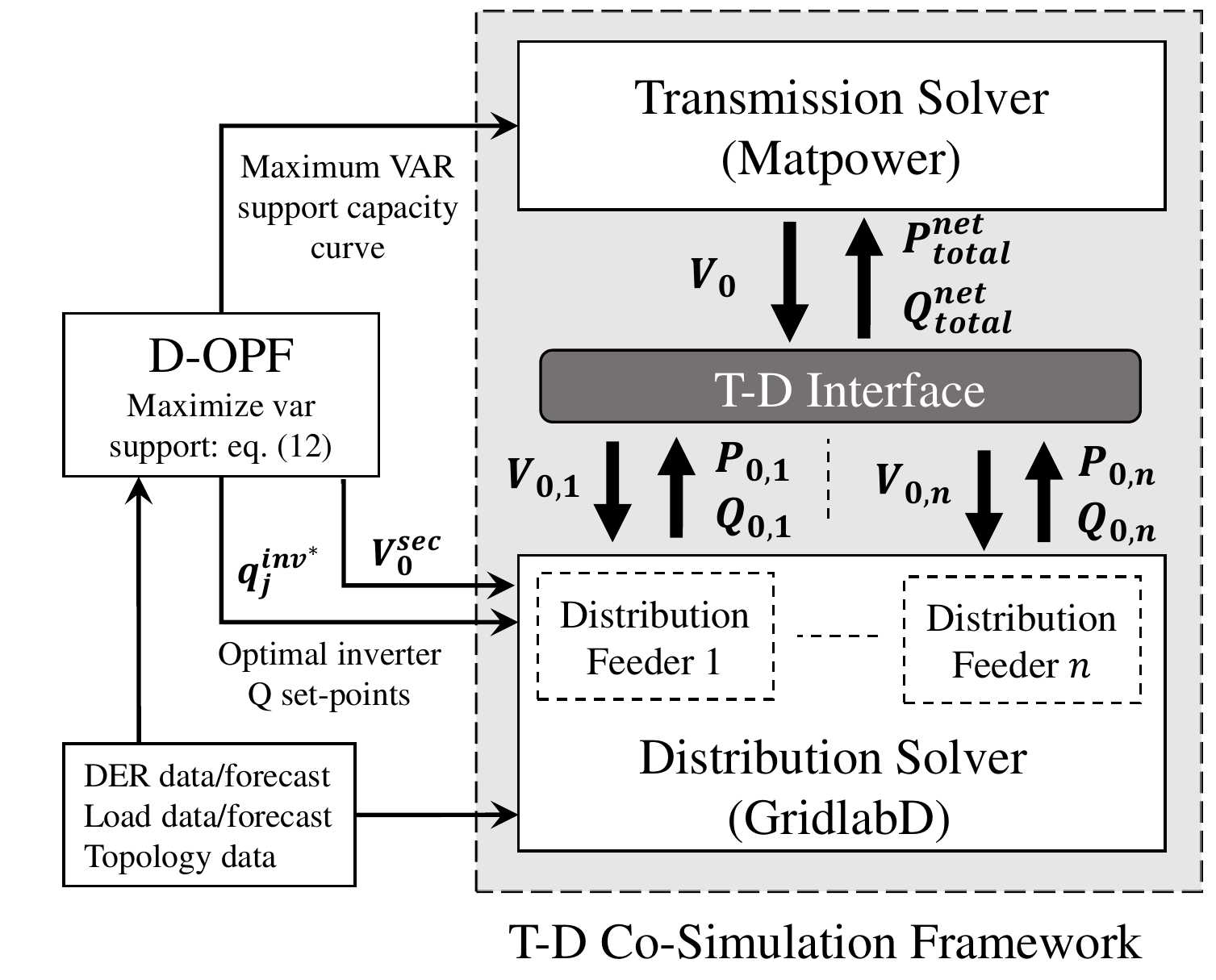}
    \vspace{-3mm}
    \caption {Proposed T-D co-simulation framework coupled with the proposed D-OPF to to assess its impact on the grid }
       \label{fig:cosimulation}
       \vspace{-4mm}
\end{figure}

A T-D co-simulation platform is developed to accurately assess the impact of the proposed var support strategy, and coupled with the D-OPF module as shown in \figurename \ref{fig:cosimulation}. The master-slave splitting (MSS) method based power flow algorithm has been used in developing this platform \cite{sun_master_2015}. We extend the MSS method to develop co-simulation for widely accepted full-scale open-source solvers i.e. Matpower for transmission and GridlabD for distribution systems. The developed platforms allows to include highly detailed unbalanced 3 phase distribution system model. The details of the method can be found in \cite{sun_master_2015}. 
Basically, it T-D interface allows the exchange of the substation voltage ($V_0$) and the net powers ($P^{net}_{total}, Q^{net}_{total}$) until convergence. Once the convergence is achieved, it interacts with the D-OPF module as shown in \figurename \ref{fig:cosimulation} i.e. receives the next set of optimal inverter var set-points and optimal secondary voltage set-point from D-OPF.  
\section{Results and Discussion}
\begin{figure}
	\centering
	\includegraphics[trim=-0in 0in 0in 0in,width=2.5in]{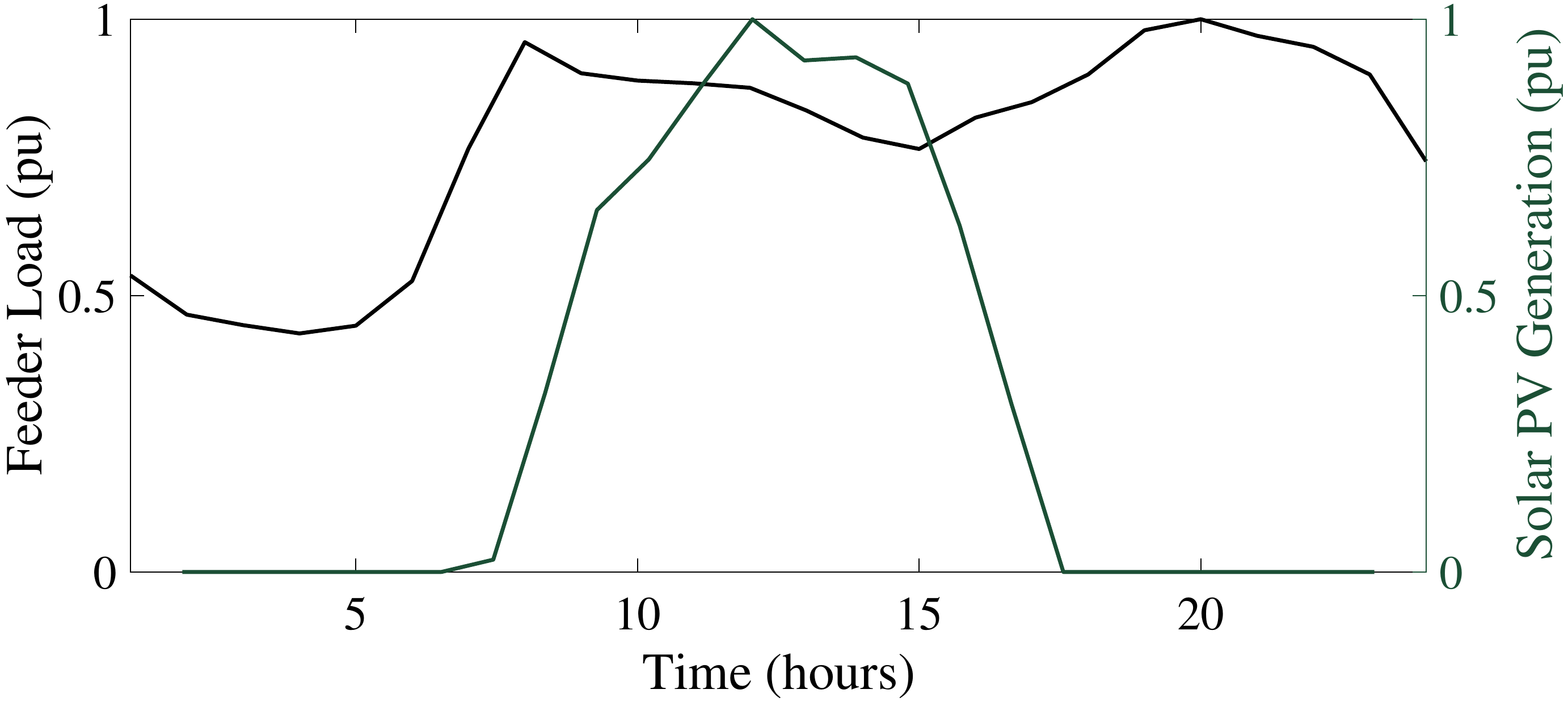}
    \vspace{-3mm}
    \caption {Normalized Daily load profile and solar PV generation profile for 24 hours with maximum value as 1 pu}
    \label{fig:daily_profile}
    \vspace{-4mm}
\end{figure}

An integrated test system is constructed by replacing 90 MW of the aggregated load at bus 5 of the 9 bus transmission system (bus T5) with several IEEE 13 bus distribution feeders to match the load. Daily load and solar PV generation curves are are normalized to their peak value as 1 pu and shown in \figurename \ref{fig:daily_profile}. 
DERs (solar PVs) are randomly distributed throughout the distribution feeder nodes. We consider several DER penetration levels from 20\% to 100\% for two cases i.e. without any solar curtailment and with 40\% solar generation curtailment. Here, we define the penetration level is a ratio of peak solar generation to peak load demand. The impact of maximum var support is verified via co-simulation.   
\subsection{Maximum Var Support Estimation }
\subsubsection{Maximum Var Support Curve and the Region} 
\begin{figure}
	\centering
    \vspace{-2mm}
	\includegraphics[trim=-0in 0in 0in 0in,width=3.5in]{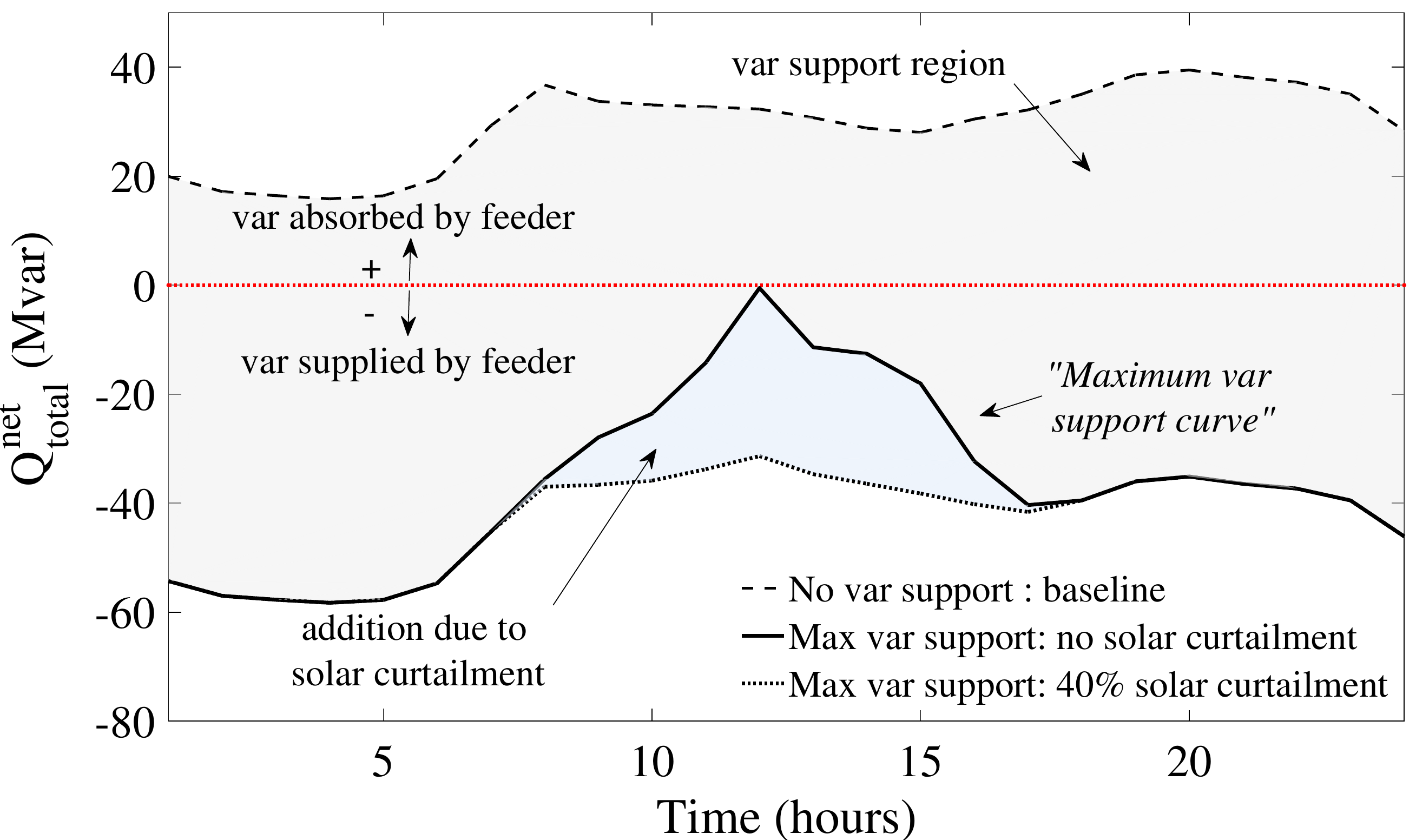}
    \vspace{-3mm}
    \caption {maximum var support region obtained through D-OPF with and without solar generation curtailment}
    \label{fig:support_region}
    \vspace{-4.5mm}
\end{figure}

\begin{figure}
	\centering
	\includegraphics[trim=-0in 0in 0in 0in,width=3.1in]{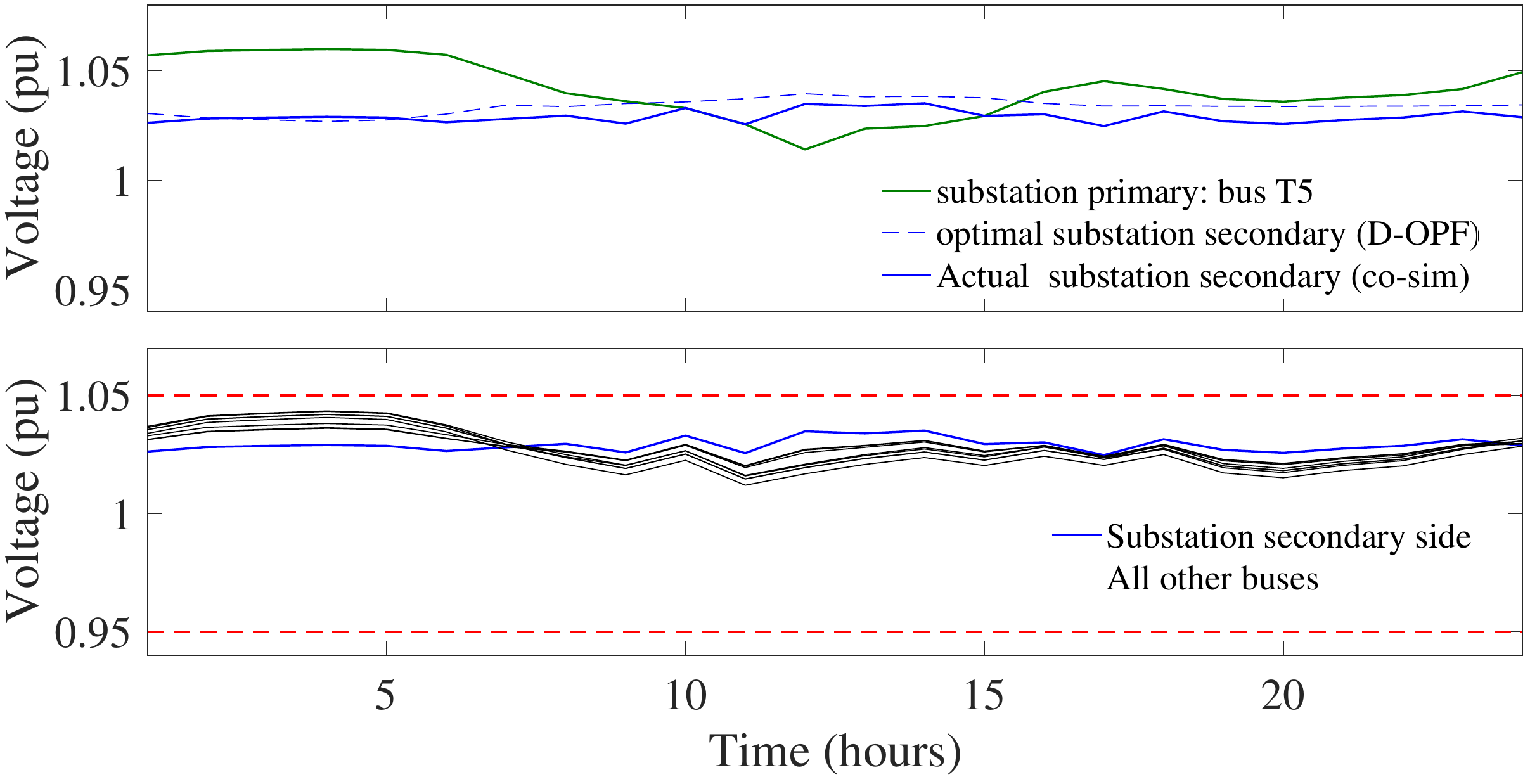}
    \vspace{-3mm}
    \caption {Distribution feeder voltages throughout the day while providing maximum var support to the grid}
    \label{fig:distvoltages_DG80}
\end{figure}

The resulting optimal var set-points and optimal secondary substation voltages from D-OPF are dispatched in co-simulation environment to obtain the true available support from all the feeders at bus T5 which is plotted as \textit{"maximum var support curve"} in \figurename \ref{fig:support_region} for 80\% DER penetration. 
Essentially, at a given time instant, DERs can meet any var demand between the dashed and the solid black curves (no var and max var support, respectively). Graphically, we term this gray shaded ares as "var support region". Note that at the peak solar time (12 noon), the var support capability is minimum as the most of the inverter capacity is occupied with the real power generation. Nonetheless, the new maximum var support curve is plotted as a dotted black line with 40\% curtailment in solar generation. This leads to the expansion in the var support region because of the addition of a new blue shaded are during solar peak. Although the solar curtailment is not economical to customers, this option exhibits the flexibility of the system. Utilizing this flexibility involves a greater discussion on policy, customer comfort, and related cost-benefit analysis. The most important aspect of this D-OPF framework is that it ensures acceptable voltage profiles in the feeders while providing maximum var support to the grid. \figurename \ref{fig:distvoltages_DG80} (top) shows the voltage profiles of primary and secondary sides of substation in co-simulation environment. It can be seen that the actual secondary side voltage (verified via co-simulation) closely follows the optimal set-points from D-OPF by changing the OLTC taps. 
\figurename \ref{fig:distvoltages_DG80} (bottom) shows all the voltages at phase A at all the feeder nodes are within the desired range. This shows the capability of the distribution system to provide promised maximum var support at any given point of time without violating its own operational limits.


\subsubsection{Impact of DER penetration levels}
\begin{figure}
	\centering
    \vspace{-2mm}
	\includegraphics[trim=-0in 0in 0in 0in,width=2.8in]{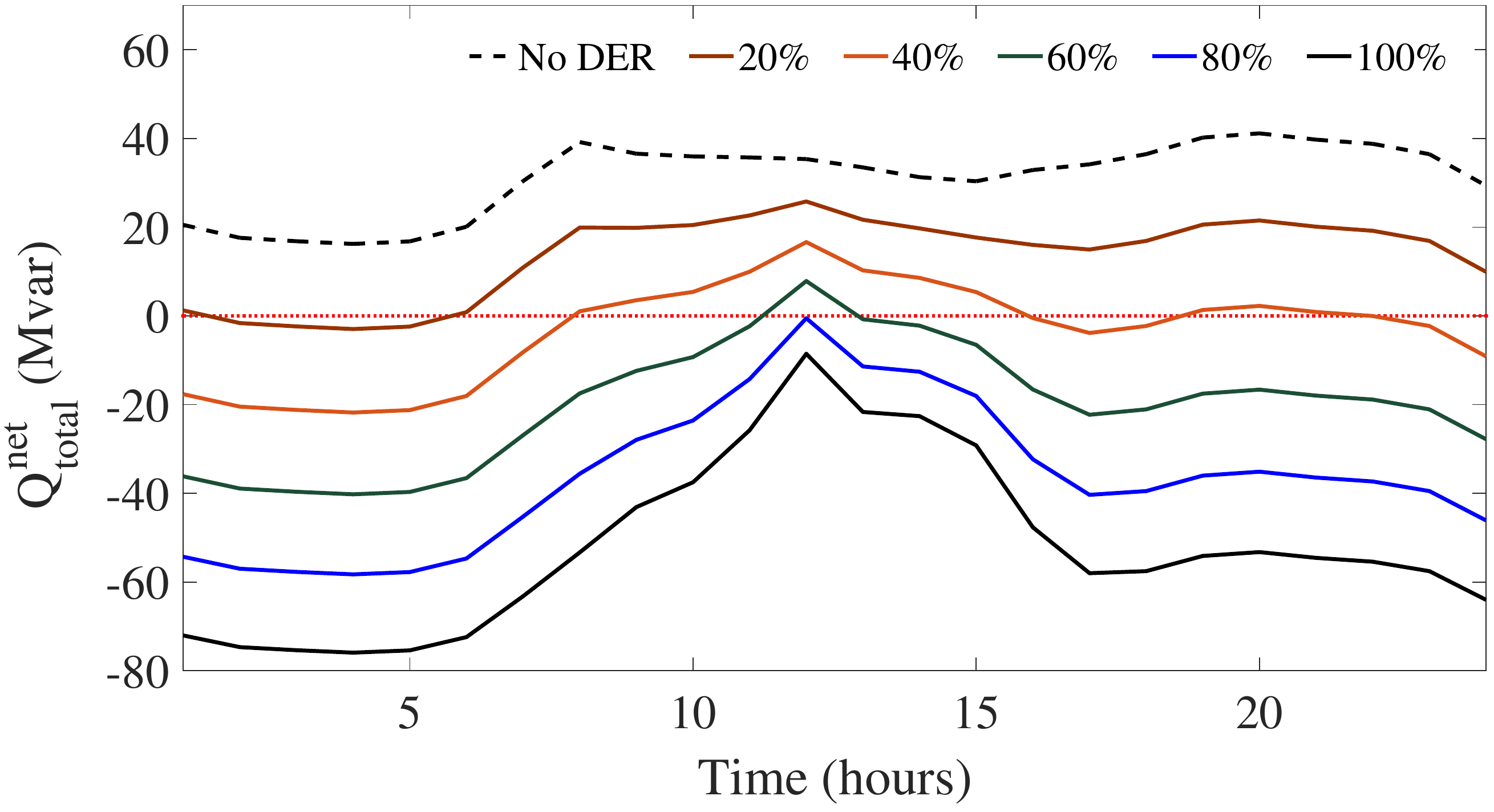}
    \vspace{-3mm}
    \caption {maximum var support curve for different DER penetration levels}
    \label{fig:varsupport_DGpen}
    \vspace{-3mm}
\end{figure}

\figurename \ref{fig:varsupport_DGpen} shows the maximum var support curves for different DER penetration levels. As the penetration level increases, the maximum var support curve goes further negative indicating there is more var available to supply by the feeder. The average \% reductions in the net feeder var demand from the no var support case are collected for various cases in the Table \ref{tab:reduction_compare}. More than 100\% reduction means the feeder can supply the var. Observe the low penetration (20\%) might not be able to supply the var to the grid, however, it still reduces the average var demand significantly (62\%) at the substation. Also, the average reduction during the solar peak (11AM-3PM) is much less than the average reduction throughout the day, in fact, it is minimum in the day because of inverter capacity is not available, however, that can be increased considerably (e.g. from 135\% to 206\% reduction with 80\% DER) by opting for solar energy curtailment.    

\begin{table}
\vspace{-4mm}
\caption{\% reduction in net feeder var demand at substation with and without solar curtailment at various DER penetration levels}
\vspace{-1mm}
\label{tab:reduction_compare}
\centering
\renewcommand{\arraystretch}{1.1}
\begin{threeparttable}
\begin{tabular}{  c  c  c | c  c  }
\hline
DER	& \multicolumn{2}{c}{No solar curtailment} & \multicolumn{2}{c}{40\% solar curtailment}  \\ \cline{2-5}

\shortstack{Penetration}  &Average & peak solar 
& Average & peak solar \\ \hline
	20\% & 62.6 & 35.52 
    & 67.23 & 53.22 \\
	40\% & 123.8 & 69.86
    & 133.1 & 105.25 \\
	60\% & 183.7 & 103.02 
    & 197.63 & 156.13 \\
	80\% & 242.3 & 135.03 
    & 260.81 & 205.86 \\
	100\% & 299.52 & 165.9 
    & 322.6 & 254.4 \\
    \hline
\end{tabular}
  \end{threeparttable}
  \vspace{-4mm}
\end{table}

\subsection{Impact of Var Support Framework on the Grid}
The proposed framework sends the maximum var support curve to the transmission grid. However, the grid might not want the maximum var support all the time; rather it can ask for the support in specific needs such as in case of voltage dips due to line contingency or to increase the load margin. We will observe two such cases with the help of the proposed co-simulation platform at the peak load scenario (8 PM).
\subsubsection{Var Support During Line Contingency}
\begin{figure}
	\centering
	\includegraphics[trim=-0in 0in 0in 0in,width=3.4in]{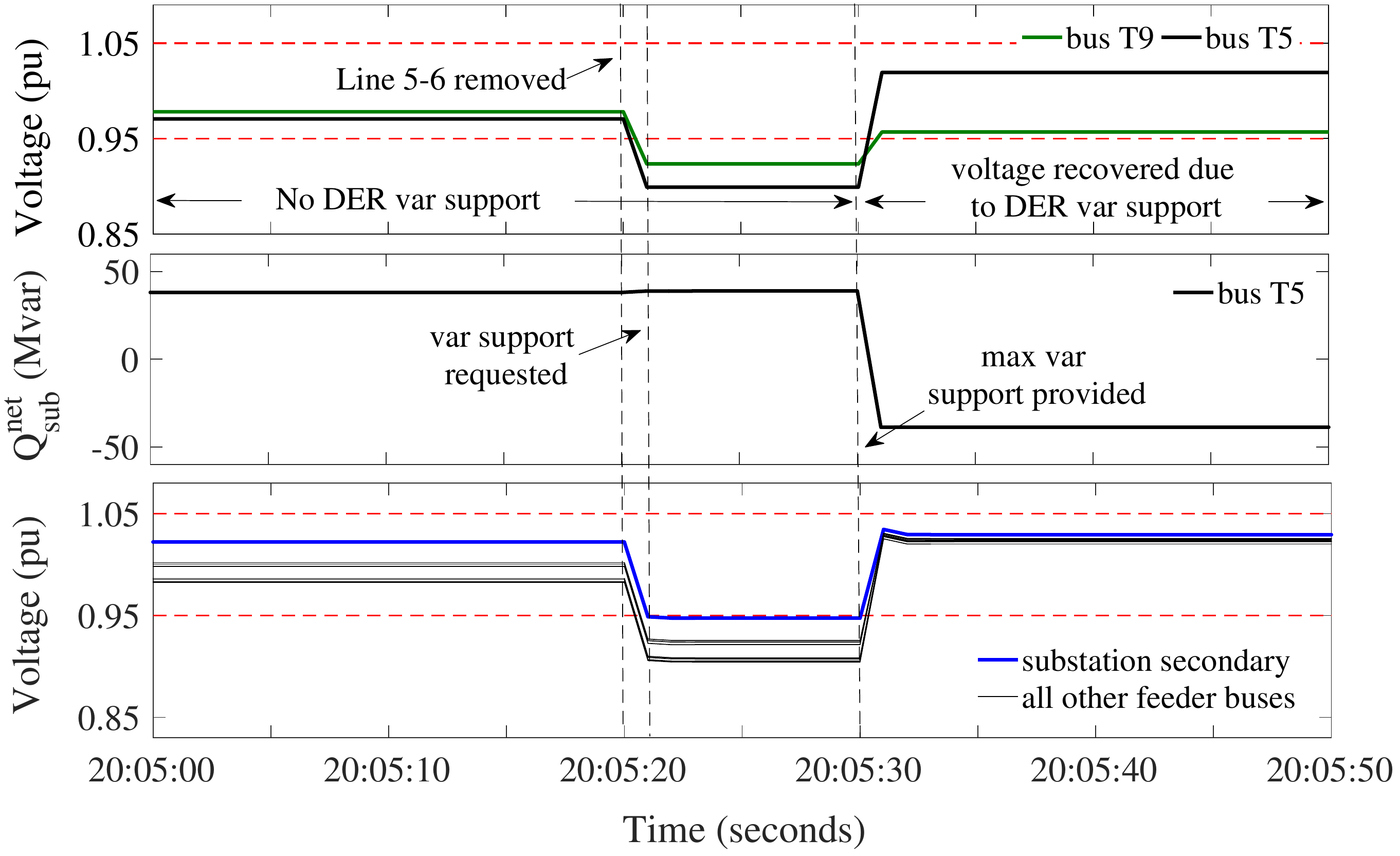}
    \vspace{-3mm}
    \caption {Var support provided by DERs in case of transmission line 5-6 contingency  }
    \label{fig:contingency56}
    \vspace{-3mm}
\end{figure}

\figurename \ref{fig:contingency56} shows how DERs var support can assist the grid in contingency situation during the peak load. The top, middle and the bottom plots show the behavior of transmission system voltages, substation net var demand, and distribution feeder voltages, respectively. Initially, the transmission system is not utilizing the var support from the DERs and operating within the limits. At point A, the transmission line 5-6 is removed which causes a sudden dip in the transmission voltages below 0.95 at bus T5 and T9. Immediately, at point B, the grid requests the var support from the DERs. Here we assume that the grid is requesting the specific var from the feeder using another grid-level OPF and the feeder can meet the request, although we do not discuss that process in this work. After a small delay of 10 seconds (communication delay, D-OPF dispatch time etc), the var support is provided at point C as can be seen in the middle plot. Note that the voltage at bus T5 gets improved more than the voltage at bus T9. This is because of the local effect of the var support provided only at T5. Simultaneously on the distribution side, the feeder voltages also experience a momentary dip as reflection of the voltage at T5. However, once OLTC operates, the voltages are brought back within the limit. This case study signifies the importance of DER var support framework as an ancillary service for the grid.

\subsubsection{Impact on the Load Margin}
\begin{figure}
	\centering
	\includegraphics[trim=-0in 0in 0in 0in,width=3in]{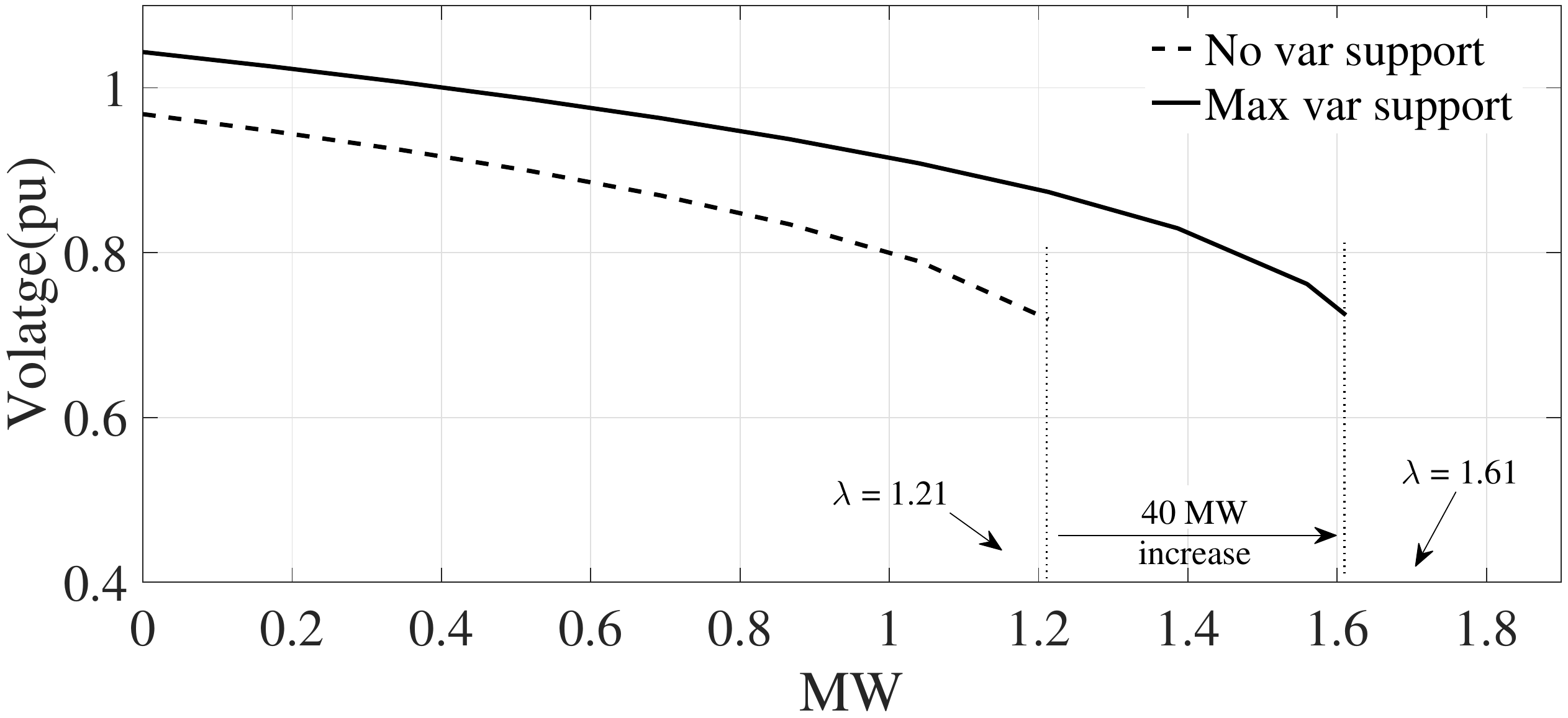}
    \vspace{-3mm}
    \caption {Comparison between $\lambda - V$ curve with no var support and maximum var support}
    \label{fig:load_margin}
    \vspace{-3mm}
\end{figure}

The var support from DERs can also help in increasing the load margin of the system from the voltage collapse point. To estimate the load margin of the integrated TD test case, the $\lambda-V$ curve (also known as PV curve) is plotted for bus T5 using co-simulation. Loads at all the distribution load buses connected to T5 are increased such that the $\lambda =0$ is the base case (90 MW) and $\lambda=1$ denotes the increase of 100 MW. Since a 3-phase unbalanced distribution feeder is used using GridlabD in the integrated system, we stop tracing the curve when the GridlabD stops converging. The PV curves with and without var support (80\% DER) are plotted in \figurename \ref{fig:load_margin}. The maximum var support increases the maximum $\lambda$ from 1.21 to 1.61 which essentially denotes the improvement of 40 MW in the load margin. 
\section{Conclusion}
In this work, we developed a DER var support framework for an integrated T-D system to test our hypothesis that the thousands of inverter-based DER devices can be exploited as geographically distributed var resources (mini-SVCs) to improve the performance of the transmission grid. There were two main questions that needed to be addressed to explore this idea. First, we proposed a D-OPF which can estimate the maximum DER var support capacity curve by while ensuring the feeder voltages are within their operational limits. This support curve then is sent to the grid to enable them take an appropriate control decision. Second, we verify the performance of the proposed framework and investigate the true impacts of var support on the grid using a co-simulation environment. The findings on an integrated T-D test system (IEEE 9 bus+IEEE 13 bus 3 phase test feeders) are encouraging and the results confirm our premise that the DERs have potential to provide volt/var ancillary services to the grid for performance enhancement, if controlled properly.

It is worth nothing that we provide a general framework for var support and impact assessment which can be extended for various other grid applications such as enabling feeder to meet the specific var request from the grid, optimizing solar curtailment, coupling with the local volt/var control at distribution side, enabling DER var support at more than one transmission buses etc. All these functionality will be explored in the future studies. 

\vspace{-0.5mm}
\bibliographystyle{IEEEtran}
\bibliography{IEEEabrv,ref_short}

\end{document}